\documentclass[11pt]{amsart}
\usepackage{amscd,amssymb,palatino}
\usepackage[utf8]{inputenc}
\usepackage[english]{babel}
\usepackage{caption}
\usepackage{etex}
\usepackage[leqno]{amsmath}
\usepackage{amssymb}
\usepackage{color}
\usepackage{shadow}
\usepackage{epsfig}
\usepackage{epic}
\usepackage{eepic}
\usepackage{graphics}
\usepackage{graphicx}
\usepackage{psfrag}
\usepackage{calc}
\usepackage{tikz}
\usepackage[dvipsnames,prologue,table]{pstricks}
\usepackage{pstcol}
\usetikzlibrary{arrows,decorations,patterns,positioning,automata,shadows,fit,shapes,calc,decorations.markings,backgrounds,scopes,decorations.text}

\usepackage{tipa}

\usepackage{fancyhdr}
\fancypagestyle{plain}{\fancyhf{}
\fancyfoot[C]{{\large\sf\bfseries\thepage}}}
\usepackage{calc}
\newcommand\1{\lower 9pt\hbox{\underbar{}}}
\numberwithin{equation}{section}

\newtheorem {Theorem}                   {Theorem}

\newtheorem {Corollary}                 {Corollary}

\theoremstyle{definition}

\newtheorem {Remark}[equation]          {Remark}

\newcommand{\pr} {\smallskip\noindent{\bf Proof\,\,}}

\begin{document}

\title{Darboux, Moser and Weinstein theorems for prequantum systems}

\author{Eva Miranda}
\thanks{{Both authors are partially supported by the Spanish State Research Agency grant PID2019-103849GB-I00 of MICIU/AEI / 10.13039/501100011033.  E. Miranda is supported by the Catalan Institution for Research and Advanced Studies via an ICREA Academia Prize 2021 and by the Alexander Von Humboldt foundation via a Friedrich Wilhelm Bessel Research Award. E. Miranda is also supported by the Spanish State
Research Agency, through the Severo Ochoa and Mar\'{\i}a de Maeztu Program for Centers and Units
of Excellence in R\&D (project CEX2020-001084-M) and partially supported by and by the AGAUR project 2021 SGR 00603. }}
\address{{Department of Mathematics} at Universitat Polit\`{e}cnica de Catalunya and Centre de Recerca Matem\`{a}tica-CRM}
\email{eva.miranda@upc.edu}
\author{Jonathan Weitsman}
\thanks{J. Weitsman was supported in part by a Simons Collaboration Grant \# 579801.}

\address{Department of Mathematics, Northeastern University, Boston, MA 02115}
\email{j.weitsman@neu.edu}
\thanks{\today}

\begin{abstract} We establish analogs of the Darboux, Moser and Weinstein theorems for prequantum systems. We show that two prequantum systems on a manifold with vanishing first cohomology, with symplectic forms defining the same cohomology class and homotopic to each other within that class, differ only by a symplectomorphism and a gauge transformation. As an application, we show that the Bohr-Sommerfeld quantization of prequantum system on a manifold with trivial first cohomology is independent of the choice of the connection.

\end{abstract}

\maketitle

\section{Introduction}

The Darboux theorem establishes that there are no local invariants in symplectic geometry. 
Namely, let $M$ be a smooth manifold and let $\omega, \omega^{\prime}$ be two symplectic forms on $M$, then:

\begin{Theorem}[Darboux \cite{darboux, darboux2}]
     {For every point $p\in M$ there exists a neighbourhood $U$ of $p$ and an embedding $\Phi:U \rightarrow M$ isotopic to the inclusion and fixing $p$ such that

$$
\Phi^{*} \omega^{\prime}=\omega|_U \text {. }
$$}
\end{Theorem}

The global invariants of a symplectic manifold include the cohomology class defined by its symplectic form. This gives a complete classification of symplectic $2$-manifolds. In general, there is the following theorem due to Moser.

\begin{Theorem}[Moser \cite{moser}]
Let $M$ be a compact manifold endowed with two symplectic forms $\omega$ and $\omega^{\prime}$. Assume that $[\omega]=\left[\omega^{\prime}\right]$,  and that there exists a path $\omega_t$ of symplectic forms such that  $[\omega_t ]= [\omega_0]$ for all $t,$ and with $\omega_0=\omega$ and 
$\omega_1=\omega^{\prime}$. Then,  there exists a diffeomorphism isotopic to the identity
 $\Phi:M \rightarrow M$ such that $\Phi^{*} \omega^{\prime}=\omega$.   
\end{Theorem}

This proof based on Moser's method can be adapted when symmetries are present.
Under the conditions above when a compact group acts on $M$ preserving $\omega$ and $\omega^{\prime},$ and where the path $\omega_t$ also consists of invariant forms, we have the following result.

\begin{Theorem}[Weinstein \cite{weinstein}]
    There exists an equivariant diffeomorphism $\Phi:M \rightarrow M$ such that $\Phi^{*} \omega^{\prime}=\omega$.
\end{Theorem}

Recall that if $(M, \omega)$ is a symplectic manifold, $\pi: L \rightarrow M$ is a (complex) line bundle, $\nabla$ is a connection on $L$, and $\operatorname{curv}(\nabla)=\omega$, the quadruple $(M, \omega, L, \nabla)$ is a \emph{prequantum system}.

The main purpose of this article is to generalize these three theorems to prequantum systems.

\begin{Theorem} \label{maintheorem1}
    
 Suppose $(M, \omega, L, \nabla)$ and $(M, \omega^{\prime}, L, \nabla^{\prime})$
are prequantum systems. Then for every point $p\in M$ there exists a neighbourhood $U$ of $p,$ an embedding $\Phi:U \rightarrow M$ isotopic to the inclusion and fixing $p,$ and a function $\phi \in C^{\infty}(U, \mathbb C)$ such that\footnote{{For the definition of the pullback of a connection under a smooth map, see e.g. \cite{ms},  p. 292, Lemma 3.}}

$$
\begin{aligned}
& \omega|_U=\Phi^{*} \omega^{\prime} \quad (\Phi^{*} L) = L|_U \\
& \nabla s=(\Phi^{*} \nabla^{\prime})s+  d \phi\otimes s, \, \text{for any} \, s\in \Gamma(L|_U).\end{aligned}
$$

If  $L$ is a hermitian line bundle and $\nabla$ is a hermitian connection, we may take  $\phi \in \mathcal{C}^{\infty}(M, i\mathbb{R}) $.
\end{Theorem}

For compact manifolds, we have

\begin{Theorem} \label{maintheorem2}
    
 Suppose $(M, \omega, L, \nabla)$ and $(M, \omega^{\prime}, L, \nabla^{\prime})$
are prequantum systems with $M$ a compact manifold. Assume that $[\omega]=\left[{\omega}^{\prime}\right]$,  and that there exists a path $\omega_t$ of symplectic forms with $[\omega_t] = [\omega_0]$ for all $t,$ such that $\omega_0=\omega$ and 
$\omega_1=\omega^{\prime}$. Assume also that $H^1(M,{\mathbb R})=0.$

Then  there exists a  a diffeomorphism  $\Phi: M\rightarrow M$ isotopic to the identity   and a function $\phi \in \mathcal{C}^{\infty}(M ,\mathbb{C})$ such that

$$
\begin{aligned}
& \omega=\Phi^{*} \omega^{\prime} \quad \Phi^{*} L = L \\
& \nabla s=(\Phi^{*} \nabla^{\prime})s+  d \phi\otimes s, \, \text{for any} \, s\in \Gamma(L).
\end{aligned}
$$

If  $L$ is a hermitian line bundle and $\nabla$ is a hermitian connection, we may take  $\phi \in \mathcal{C}^{\infty}(M, i\mathbb{R}) $.
\end{Theorem}

When there are additional symmetries the following result holds.

\begin{Theorem}\label{maintheorem3} Suppose $G$ is a compact Lie group, and suppose $M$ is a compact $G$-space. Suppose that $(M, \omega, L, \nabla)$ and $(M, \omega^{\prime}, L, \nabla^{\prime})$ are $G$-invariant prequantum systems. Assume that $[\omega]=\left[\omega^{\prime}\right]$,  and that there exists a path $\omega_t$ of $G$-invariant symplectic forms with $[\omega_t] = [\omega_0]$ for all $t,$  such that $\omega_0=\omega$ and 
$\omega_1=\omega^{\prime}.$ Assume also that $H^1(M,{\mathbb R})=0$.

Then  there exists a   $G$-equivariant diffeomorphism  $\Phi: M\rightarrow M$ isotopic to the identity   and a $G$-invariant function $\phi \in \mathcal{C}^{\infty}(M ,\mathbb{C})$ such that

$$
\begin{aligned}
& \omega=\Phi^{*} \omega^{\prime} \quad \Phi^{*} L = L \\
& \nabla s=(\Phi^{*} \nabla^{\prime})s+  d \phi\otimes s, \, \text{for any} \, s\in \Gamma(L).
\end{aligned}
$$

If  $L$ is a hermitian line bundle and $\nabla$ is a hermitian connection, we may take  $\phi \in \mathcal{C}^{\infty}(M, i\mathbb{R}) $.
\end{Theorem}

\begin{Remark}
  The function $\phi$ is called a \emph{gauge transformation}. We have thus proven that two prequantum systems, with associated symplectic forms lying in the same cohomology class and homotopic  to each other, differ only by a symplectomorphism and a gauge transformation.  Observe that it is not possible to dispense with the gauge transformation:  Suppose $\nabla$ is the trivial connection on the trivialized bundle $L = M \times {\mathbb C},$ and $\nabla^\prime$ is some other connection on $L.$ Then for any section $s \in \Gamma(M,L),$ we have $\nabla^\prime s = \nabla s +d\phi\otimes s,$ where $\phi$ is somewhere non-constant.   And no symplectomorphism will transform $\phi$ into an everywhere constant function (since diffeomorphisms take constant functions to constant functions).
\end{Remark}
 
As an application of this theorem, we prove that for  manifolds with trivial first cohomology, the Bohr-Sommerfeld quantization of a prequantum system $(M,\omega,L,\nabla)$ does not depend on the choice of the connection $\nabla;$  in other words if $(M,\omega,L,\nabla^\prime)$ is another prequantum system associated with the same manifold, line bundle, and symplectic form, its Bohr-Sommerfeld quantization coincides with that of $(M,\omega,L,\nabla).$
 If $H^{1}(M) \neq 0$ a similar result holds if the holonomy representations of $\nabla$ and $\nabla^{\prime}$ coincide.
 
 {\bf Acknowledgement:}  We would like to thank Peter Crooks for helpful discussions on this paper.

\section{Proof of the main theorems}
We first prove Theorem \ref{maintheorem1}.
\begin{proof}

 By the Darboux theorem, we know that there exists a neighborhood $U$ of $p,$ which we may take to be diffeomorphic to an open ball, and {an embedding $\phi: U \rightarrow M$
isotopic to the inclusion preserving $p$ such that 
$$
{\Phi^{^*} \omega^{\prime}=\omega|_U .}
$$
Since} $\Phi$ is isotopic to the inclusion we also have

$$
\Phi^{*} L=L|_U
$$

It remains to compare the two connections and consider the difference
{
$$
\nabla s-(\Phi^{*} \nabla^{\prime}) s \, \text{for any} \, s\in\Gamma(L|_U).
$$}

Since the space of connections on a line bundle on $U$ is an affine space modeled on $\Omega^{1}\left(U, \mathbb C\right) ,$ we know that

$$
\nabla s-(\Phi^{*} \nabla^{\prime})s= \xi\otimes s
$$

for some $\xi \in \Omega^{1}\left(U, \mathbb C\right)$.

Since  $\operatorname{curv}\nabla=\operatorname{curv}\Phi^{*} \left(\nabla^{\prime}\right)$, it follows
that $d \xi=0$. Since $H^{1}(U,{\mathbb R})=0$,
$
\xi=d \phi$ for  some $\phi \in \mathcal{C}^{\infty}(U,\mathbb{C}).
$
If the line bundle $L$ is hermitian and the connection is unitary, we may repeat this argument to show that we may take $\phi \in \mathcal{C}^{\infty}(U,i\mathbb{R}).$
\end{proof}

We next prove Theorem \ref{maintheorem2}.

\begin{proof} The existence of the diffeomorphism is guaranteed by Moser's theorem. 
We again have $\Phi^{*} L=L$, and again there exists $\xi \in \Omega^{1}\left(M, \mathbb C\right)$ so that
$$
\nabla s-(\Phi^{*} \nabla^{\prime})s= \xi\otimes s \, \text{for any} \, s\in\Gamma(L).
$$
Again $d \xi=0$, so since $H^1(M,{\mathbb R}) = 0,$ $\xi=d \phi$ for some $\phi,$ as needed. If the line bundle $L$ is hermitian and the connection is unitary, we may repeat the argument to show that we may take $\phi \in \mathcal{C}^{\infty}(M,i\mathbb{R}).$
    
\end{proof}

We now proceed to the proof of Theorem \ref{maintheorem3}.

\begin{proof} The existence of a $G$-equivariant diffeomorphism $\Phi$ is given by Weinstein's theorem.

Since $H^1(M,{\mathbb R})=0,$  again there exists a function $\phi$ such that

\begin{equation}\label{eqn1}
 \nabla s-(\Phi^{*} \nabla^{\prime})s= d\phi\otimes  s  ~  \text{for any} \, s\in\Gamma(L).
\end{equation}

Let $\bar{\phi}$ be the average of $\phi$ over $G$ using the Haar measure $\mu$, i.e., $$\bar{\phi}:=\int_G \alpha_g^*\phi~ d\mu(g),$$ where we have denoted by $\alpha_g$ the diffeomorphism giving the action of $g \in G$ on $M.$

Then since $\nabla$ and $\nabla^{\prime}$ are $G$-invariant, we obtain

$$
 \nabla s-(\Phi^{*} \nabla^{\prime})s=  d \bar{\phi}\otimes s 
$$

by averaging both sides of equation (\ref{eqn1}).

Again, if the line bundle $L$ is hermitian and the connection is unitary, we may repeat this argument to show that we may take $\phi \in \mathcal{C}^{\infty}(M,i\mathbb{R})^G.$
    
\end{proof}
\begin{Remark}\label{remark}(The case where $H^1(M,{\mathbb R}) \neq 0$)
Suppose we are given two connections   $\nabla$ and $\nabla^{\prime}$ on $L$ with the same curvature. Then for any section $s \in \Gamma(M,L),$ we have $\nabla s=\nabla^{\prime} s +  \alpha \otimes s$ for some $\alpha \in \Omega^1(M, {\mathbb C})$ with $d \alpha = 0$.

If $H^1(M,{\mathbb C}) \neq 0,$ we cannot conclude that $\alpha = d\phi$ for some function $\phi.$ However, if  $\nabla$ and $\nabla^{\prime}$  also have the same holonomy representation, it follows that the holonomy representation of $\alpha,$ thought as a connection on the trivialized bundle $L\otimes L^*=\mathbb C$, is trivial.  
That means that for every oriented curve $C$ in $M$, $\int_C \alpha = 2 \pi n$, where $n \in \mathbb Z$.  This does not imply that $\alpha$ is exact; but there exists a map $t: M \to S^1$ with $\alpha = t^*( d \theta)$, where $d\theta$ is the translation invariant one form on $S^1$.\footnote{Morally $\alpha = t^{-1} dt,$ where $t : M \to S^1$ is a (possibly homotopically nontrivial) element of the gauge group ${\rm Map}(M,S^1).$}
    
\end{Remark}

\section{Application to Quantization in a Real Polarization}

Let $(M, \omega, L, \nabla)$ be a prequantum system. Recall that a Lagrangian submanifold $\Lambda$
is {\em integral} if $\left.(L, \nabla)\right|_{\Lambda}$ is trivial as a bundle with connection.

\begin{Corollary}\label{cor1} Let  $M$ be a  compact manifold with $H^1(M,{\mathbb R}) = 0$ and let  $(M, \omega, L, \nabla)$ and $(M, \omega, L, \nabla^{\prime})$ be two prequantum systems. 
A Lagrangian submanifold $\Lambda\subset M$ is integral for $(M, \omega, L, \nabla)$ if and only if it is integral for $(M, \omega, L, \nabla^{\prime})$.

\end{Corollary}

In particular, suppose we are given a real polarization of $M$, that is, a foliation by Lagrangian submanifolds.  Consider from now on prequantum systems $(M, \omega, L, \nabla)$  where where $L$ is hermitian and the connection $\nabla$ is unitary.
The Bohr-Sommerfeld quantization of  $(M, \omega, L, \nabla)$ associated to this polarization is given by the vector space generated by the integral leaves of the foliation, or equivalently by the covariant constant sections of   $(L,\nabla)$ on those leaves.

{In  \cite{guilleminsternberg} Guillemin and Sternberg consider the case where the foliation is given by a map $\pi : M \to B,$ where $B$ is a simply connected subset of $ {\mathbb R}^n$ and the map $F: M\longrightarrow B$ is the moment map for an integrable system away from the singularities of the foliation.  The Bohr-Sommerfeld quantization is then determined by the integral points of the image of the moment map $F$.}

Then Corollary 1 implies

 \begin{Corollary}\label{cor2}
   Let $(M, \omega)$ be a symplectic manifold with $H^1(M,{\mathbb R})=0.$  Let $L$ be a complex line bundle on $M$ with $c_1(L) = [\omega].$  Choose a connection $\nabla$ on $L$ so that $(M,\omega,L,\nabla)$ is  a prequantum system. Suppose we are given a real polarization of $M.$ Then the quantization of   $(M, \omega, L, \nabla)$ with this real polarization is independent of the choice of $\nabla.$
 \end{Corollary}
 
 \begin{Remark}  Note that in the case where the foliation is given as in  \cite{guilleminsternberg} by an integrable system, the generic leaves are tori, on which $(L,\nabla)$ may have complicated holonomy representations.  We have shown that the global condition $H^1(M,{\mathbb R}) = 0 $ guarantees these holonomy representations are independent of $\nabla.$\end{Remark}

\begin{Remark} For complex polarizations, an analog of Corollary \ref{cor2} holds, since the index of the Dolbeault $\bar\partial$ operator on $(M, L)$ depends only on $M$ and $c_{1}(L)=[\omega]$, not on the connection; this is due to the Riemann-Roch theorem

\begin{equation}\label{eqn:index}
  \operatorname{ind}(\bar\partial)=\int_{M} Td(T M) e^{[\omega]}.   
\end{equation}\end{Remark}

Note that the right-hand side of equation (\ref{eqn:index}) depends only on the class $[\omega]$ and not on $\nabla$.

\begin{Remark}
    In the case where $H^1(M,{\mathbb R}) \neq 0,$ any two connections $\nabla$ and $\nabla^{\prime}$ with the same curvature still give the same quantization if, in addition, the trivial bundle $L \otimes L^*$, equipped with the flat connection arising from $\nabla$ and $\nabla^{\prime}$, is trivial as a bundle with connection. 
\end{Remark}

\end{document}